\documentclass[12pt]{amsart}
\usepackage{amssymb,latexsym}
\usepackage{mathrsfs}
\usepackage{amsmath}
\usepackage{layout}
\usepackage{setspace}
\usepackage{nicefrac}
\theoremstyle{plain}

\theoremstyle{definition}

\usepackage{amsthm}  
\newtheoremstyle{named}{}{}{\itshape}{}{\bfseries}{.}{.5em}{\thmnote{#3 }#1} 
\theoremstyle{named} 
\newtheorem*{namedtheorem}{Theorem} 

\numberwithin{equation}{section}
\numberwithin{lemma}{section}
\numberwithin{theorem}{section}
\numberwithin{proposition}{section}
\begin{document}
\thispagestyle{empty}
\onehalfspacing
\title{Notes on Thin Matrix Groups}
\author{Peter Sarnak}
\address{
Deptartment of Mathematics\\
Princeton University \\
Princeton, NJ   08544-100 USA }

  \email{sarnak@math.princeton.edu}
\maketitle

\begin{abstract}
   These notes were prepared for the MSRI hot topics workshop on superstrong approximation (2012). We give a brief overview of the developments in the theory, especially the fundamental expansion theorem. Applications to diophantine problems on orbits of integer matrix groups, the affine sieve, group theory, gonality of curves and Heegaard genus of hyperbolic three manifolds, are given. We also discuss the ubiquity of thin matrix groups in various contexts, and in particular that of monodromy groups.*
\end{abstract}\footnote{There is some overlap with Lubotzky's recent Colloquium Lectures \cite{Lu1} and Green's note \cite{Gr} on group theoretic combinatorics.}
\vspace{.15in}

\section{The Fundamental Expansion Theorem}
The Chinese Remainder Theorem for $SL_n (\mathbb{Z})$ asserts, among other things, that for $q \geq 1$, the reduction $\pi_q : SL_n (\mathbb{Z}) \longrightarrow SL_n (\mathbb{Z} / q \mathbb{Z} )$ is onto. Far less elementary is the extension of this feature to $G(\mathbb{Z} )$ where $G$ is a suitable matrix algebraic group defined over $\mathbb{Q}$. The general form of this phenomenon for arithmetic groups is known as strong approximation and it is well understood \cite{P-R}.\\
\indent There is a quantification of the above that is not as well known as it should be, as it turns out to be very powerful in many contexts. We call this \textquotedblleft superstrong\textquotedblright \ approximation and it asserts that if we choose a finite symmetric ($s \in S$ \ \ if \ $s^{-1} \in S$) generating set $S$ of $SL_n (\mathbb{Z} )$, then the congruence Cayley graphs $(X_q , S)$ form an expander family as $q$ goes to infinity (see \cite{H-L-W} for the definition and properties of expanders). Here the vertices $x$ of the $| S |$-regular connected graph $(X_q , S)$ are the elements of $SL_n (\mathbb{Z} / q \mathbb{Z} )$ and the edges run from $x$ to $sx, \ s \in S$. The proof of this expansion property for $SL_2 (\mathbb{Z} )$ has its roots in Selberg's 3/16 lower bound for the first eigenvalue $\lambda_1$ of the Laplacian on the hyperbolic surface $\varGamma \backslash \mathbb{H} , \ \varGamma$ a congruence subgroup of $SL_2 (\mathbb{Z} )$, (\cite{Se}). This bound is an approximation to the Ramanujan/Selberg Conjecture for automorphic forms on $GL_2 / \mathbb{Q}$. The generalizations of the expansion property to $G(\mathbb{Z} )$ where $G$ is say a semisimple matrix group defined over $\mathbb{Q}$ is also known thanks to developments towards the general Ramanujan Conjectures that have been established (\cite{B-S}, \cite{Cl}, \cite{Sa1}). This general expansion for these $G(\mathbb{Z})$'s also goes by the name \textquoteleft property $\tau$\textquoteright \ for congruence subgroups \cite{Lu2}.\\
\indent Let $\varGamma$ be a finitely generated subgroup of $GL_n (\mathbb{Z} )$ (more generally later on we allow it to be in $GL_n (K)$ where $K$ a number field) and denote its Zariski closure: $Zcl(\varGamma )$, by $G$. If $\varGamma$ is of finite index in $G(\mathbb{Z} )$, then the discussion above of strong and superstrong approximation can be applied. However, if $\varGamma$ is of infinite index in $G(\mathbb{Z} )$, then Vol $( \varGamma \backslash G (\mathbb{R} ) ) = \infty$ \ and the techniques used to prove both of these properties don't apply. In this case we call $\varGamma \ $\textquotedblleft thin\textquotedblright . It is remarkable that under suitable natural hypotheses, strong approximation continues to hold in this thin context. The first result in this direction is \cite{M-V-W}, and Weisfeiler extended it much further. More recent and effective treatments of this can be found in \cite{N1} and \cite{L-Pi}. An example of the statement of strong approximation in this context is: suppose that $Zcl(\varGamma ) = SL_n$ , then there is a $q_0 = q_0 (\varGamma )$ such that for $(q,  q_0 ) = 1, \ \ \pi_q : \varGamma \longrightarrow SL_n (\mathbb{Z} / q \mathbb{Z})$ is onto. That the expansion property might continue to hold for thin groups was first suggested by Lubotzky in the 90's (\cite{L-W}). Thanks to a number of major developments by many people (\cite{S-X}, \cite{Ga}, \cite{He}, {\cite{B-G1}, \cite{B-G-S1}, \cite{P-Sz}, \cite{B-G-T}, \cite{V}), the general expansion property is now known. The almost final version (almost because of the restriction that $q$ be squarefree) is due to Salehi and Varju \cite{S-V}.\\
\begin{namedtheorem}[The Fundamental Expansion] 
 Let $\varGamma \leq SL_n (\mathbb{Q})$ be a finitely generated group with a symmetric generating set $S$. Then the congruence graphs $(\pi_q (\varGamma ), S)$, for $q$ squarefree and coprime to a finite set of primes (which depend on $\varGamma$ ), are an expander family,  iff $G^o$ the identity component of $G : = Zcl(\varGamma )$, is perfect (i.e. $[G^o, G^o] = G^o)$. Moreover the determination of the expansion constant is in principle effective\textemdash if not feasible.\footnote{Matching to some extent in this general setting the quality of expansion that is known when $\varGamma$ is arithmetic, remains an open problem.}
\end{namedtheorem} 
I will not review the techniques leading to the proof of this theorem (they have been discussed in many places including Kowalski and Tao's blogs) other than to point out that it involves three steps, the opening, the middlegame and the endgame. The endgame establishes the expansion by combining sufficiently strong (but still quite crude) upper bounds for the number of closed circuits in these graphs with largeness properties of the dimensions of the irreducible representations of the finite groups $G(\mathbb{Z} / q \mathbb{Z})$. In some cases (indeed all for which reasonable bounds for the expansion are known) the proof involves the endgame only (\cite{S-X}, \cite{Ga}). In the general case, the upper bounds for the number of closed circuits is derived combinatorially. The opening and middlegame involve showing that smaller subsets of $G(\mathbb{Z} / p \mathbb{Z})$ grow substantially when multiplied by themselves at least three times (\cite{He} and the extensions \cite{P-Sz} and \cite{B-G-T}). A critical ingredient in the early treatments was the \textquoteleft sum-product\textquoteright \ theorem (\cite{B-K-T}) in finite fields. The middlegame is concerned with moderately large sets and is further handled by the crucial \textquotedblleft flattening lemma\textquotedblright \ \cite{B-G1}. The latter also has its roots in combinatorics appealing to the Balog-Szemeredi Theorem (\cite{Ba-Sz}, \cite{Go}). When $q$ is not prime, the analysis and combinatorics is far more complicated and difficult due to the many subgroups of $G(\mathbb{Z} / q \mathbb{Z})$. It is handled in \cite{B-G-S1} for $SL_2$ and in \cite{V} in general. 
\section{Applications}
\subsection{The affine sieve and diophantine analysis}
The impetus for developing the expansion property for thin groups arose in connection with diophantine problems (in particular sieve problems for values of polynomials) on orbits of such thin groups (\cite{B-G-S1}). Both strong approximation and superstrong approximation are crucial ingredients in executing a Brun combinatorial sieve in this setting. The theory is by now quite advanced and in particular the basic theorem of the affine sieve has been established in all cases where it is expected to hold (\cite{S-S}).\\
\indent For various special examples, such as for integral Apollonian packings, which has turned out to be one of the gems of the theory (\cite{Sa2}), much more can be said thanks to special features. Firstly, in this case one can develop an archimedian count for the number of points in an orbit in a large region. This is done by combining spectral methods (using techniques which when $\varGamma$ is a geometrically finite subgroup of $O(n-1, 1)(\mathbb{R} )$ go back to \cite{Pa}, \cite{Su} and \cite{L-P}) with ergodic theoretic methods (\cite{K-O}, \cite{O-S}, \cite{L-O}, \cite{Vino}). For the diophantine applications, one needs an archimedian spectral gap for the induced congruence groups, rather than the combinatorial expansion. \cite{B-G-S2} establishes the transfer of this information from the combinatorial to archimedian setting in this infinite volume case.\\
\indent Two recent highlights of these developments are the \textquoteleft almost all\textquoteright \ local to global results \cite{B-K1} and \cite{B-K2}. The first concerns integral Apollonian packings and the question is which numbers are curvatures? The expected local to global conjecture (\cite{G-L-M-W-Y}, \cite{F-S}) is proven for all but a zero density set of integers (the Conjecture asserts that there are only a finite number of exceptions). Prior to that \cite{B-F} had shown that the number of integers that are achieved is of positive density. The second development concerns the Zaremba problem which asserts that if $A \geq 5$, the set of integers $q \geq 1$ for which there is a $1 \leq b \leq q, \ (b, q) = 1$ and for which the coefficients of the continued fraction of $b / q$ are bounded by $A$, consists of all of $\mathbb{N}$. In \cite{B-K2} the theory of thin subgroups of $SL_2 (\mathbb{Z} )$ is extended to thin semi (sub) groups (one has to abandon direct spectral methods and replace them by dynamical ones \cite{Lal}, \cite{B-G-S2}). In \cite{B-K2} it is shown  that for $A \geq 50$, the set of exceptions to the Zaremba conjecture is of zero density in $\mathbb{N}$.
\subsection{Random Elements in $\varGamma$}
It is well known that for any reasonable notion of randomness, the random $f \in \mathbb{Q} [x] $ is irreducible and has Galois group the full symmetric group on the degree of $f$ symbols. In \cite{Ri} the study of such questions for the characteristic polynomial $f_\gamma$ of a random element $\gamma$ in $Sp(2g, \mathbb{Z} )$ and more general $\varGamma$'s, was initiated. The random element in $Sp(2g, \mathbb{Z} )$ is generated by running a symmetric random walk with respect to a measure $\mu$ whose support generates  $Sp(2g, \mathbb{Z} )$. The expansion property is used via a sieving argument to show that the probability that $f_\gamma$ is reducible is exponentially small. This and some generalizations are then coupled with the theory of the mapping class group $M$ to show that the random element in $M$ is pseudo Anosov. These irreducibility questions and much more,\footnote{For example to showing that it is very unlikely that a random three dimensional manifold in the Dunfield-Thurston model \cite{D-T}, has a positive first Betti number.} are extended and refined especially in terms of the sieves that are applied, in the monograph \cite{Ko}. Again, strong and superstrong approximation plays a central role.\\
\indent In a different direction \cite{L-M} examine some group theoretic questions for linear groups using a random walk and a sieve. An example of what they show is: Let $\varGamma$ be a finitely generated subgroup of $GL_n (\mathbb{C})$ which is not virtually solvable, then the set of proper powers $P : = \operatorname*{\cup}_{m=2}^\infty\limits \{ \gamma^m : \gamma \in \varGamma \}$, is exponentially small (in terms of hitting $P$ in a long random walk). In particular, this resolved an open question as to whether finitely many translates of $P$ can cover $\varGamma$, the answer being no.
\subsection{Gonality and Heegaard Genus}
A compact Riemann surface of genus $g$ can be realized as a covering of the plane of degree at most $g + 1$ (Riemann-Roch). The gonality $d(X)$ of $X$ is the minimal degree of such a realization. Unlike $g(X), \ d(X) \ $ is a subtle conformal invariant. In \cite{Z} (and later \cite{A1}) the differential geometric inequality of \cite{Y-Y} is extended to the setting of $X = \varGamma \backslash \mathbb{H}$, \ a finite area quotient (orbifold) of the hyperbolic plane. If $A(X)$ is its area and $\lambda_1 (X)$ its first Laplace eigenvalue, then 
\begin{equation}
   d(X) \geq \frac{\lambda_1 (X) A(X)}{8 \pi} .
\end{equation} 
This together with the known bounds towards the Ramanujan/Selberg conjectures for congruence (arithmetic) $X$'s (see \cite{B-B} for the best bounds for $GL_2/K, \ K$ a number field which is what is relevant here) imply that for these $X$'s, the ratio of any two of $d(X), A(X)$ and $(g(X)+1)$ is bounded universally from above and below.\\

There is a generalization of (2.1) to finite volume quotients $X = \varGamma \backslash \mathbb{H}^m$ (orbifolds) of hyperbolic $m$-space \cite{A-B-S-W}. This is stated in terms of \cite{L-Y}'s notion of conformal volume. It gives an inequality between Vol$ (X), \ \lambda_1 (X)$ and the conformal volume of a piecewise conformal map of $X$ into $S^n$. Again, this together with the known universal lower bounds for $\lambda_1 (X)$ when $X$ is congruence arithmetic (\cite{B-S}, \cite{Cl}) gives a linear in the volume, lower bound for the conformal volume of a conformal map of $X$ to $S^m$. This has a nice application to reflection groups. A discrete group of motions of $\mathbb{H}^m$ is called a reflection group if it is generated by reflections (a reflection of $\mathbb{H}^m$ is a nontrivial isometry which fixes an $m-1$ dimensional hyperplane). Using the inequalities mentioned above, one shows (\cite{L-M-R} \ for $m=2$ \ and \cite{A-B-S-W} \ for $m > 2$) that the set of maximal arithmetic reflection groups is finite for each $m$. Now Vinberg \cite{Vi} and \cite{P} have shown that for $m \geq 1000$, a reflection group can never be a lattice. Thus the totality of all maximal arithmetic reflection groups is finite.\\

(2.1) has interesting applications to diophantine equations. As observed in \cite{A2} and \cite{Fr}, Faltings' finiteness theorem for rational points on subvarieties of abelian varieties \cite{Fa} can be used to prove finiteness of rational points on curves, whose coordinates lie in the union of all number fields of a bounded degree, as long as one can show the gonality of the curve is large enough. For example, if $X_0 (N) / \mathbb{Q}$ is the familiar modular curve of level $N$ and if $D$ is given, then for $N \geq 230 D$\footnote{This follows from (1) and explicit Ramanujan bounds.}, the set of points on $X_0 (N)$ with coordinates in the union of all number of fields of degree at most $D$, is finite! Recently \cite{E-H-K} have applied similar reasoning to a diophantine problem on a tower of curves. It arises from questions of reducibility and symmetry of specializations of members of a 1-parameter family of varieties. The curves that arise (as the parameter) are determined by the monodromy group $\varGamma$ of the family (see below), and it lies in $Sp(2g, \mathbb{Z} )$ and is assumed to be Zariski dense in $Sp(2g)$. In order to show that the gonalities of the curves in question increase quickly enough, they use the combinatorial expansion that is provided by the Fundamental Expansion Theorem. Typically it is not known if $\varGamma$ is thin or not in this context(see Section 3), but the beauty of the Fundamental Theorem is that one does not need to know!\\
\indent There is an inequality similar to (2.1) for the Heegaard genus of a hyperbolic 3-manifold $X$. It is known that such an $X$ can be decomposed into two handle bodies with common boundary a surface of genus $h$ (called a Heegaard splitting). The minimal genus of such a surface in a splitting is called the Heegaard genus of $X$ which we denote by $g(X)$. Like the gonality, it is a much more subtle (this time topological) invariant of $X$ than its volume. In \cite{La} (see Theorem 4.1 and [Bu]) it is shown that for complete $X$ of finite volume
\begin{equation}
   g(X) \geq \frac{\min [\lambda_1 (X), 1] \cdot \text{Vol} (X)}{32 \pi} .
\end{equation} 
Applying this together with the universal lower bounds for $\lambda_1$ for congruence arithmetic $X$'s, shows that the Heegaard genus of a congruence hyperbolic three manifold is in order of magnitude, a linear function of its volume. In particular, any arithmetic 3-manifold has an infinite tower (by congruence subgroups) of coverings whose Heegaard genus grows linearly with the volume. One can ask if the same is true for any hyperbolic 3-manifold and the answer is yes as was shown in \cite{L-L-R}. Using local rigidity of lattices in $SL_2 ( \mathbb{C} )$ one can realize $\varGamma$ where $X = \varGamma \backslash \mathbb{H}^3$, as a finitely generated subgroup of $SL_2 (K)$, where $K$ is some number field. If $\varGamma$ is not arithmetic then $\varGamma$ is thin (in $SL_2 (O_K )$ perhaps allowing denominators at finitely many places), since its projection on the identity embedding of $K $ into $ \mathbb{C}$ is discrete. Using the Fundamental Expansion Theorem gives a lower bound on $\lambda_1$ for a \textquoteleft congruence tower\textquoteright \ of $\varGamma$ and one then applies (2.2).\\
\indent A related application of the expansion is to some questions in knot theory. Answering a question of Gromov, Pardon \cite{Par} recently showed that there are isotopy classes of knots in $S^3$ which have arbitrary large distortion. In fact he shows that torus knots have this property. In \cite{Gr-Gu} a large family of knots with large distortion is constructed using hyperbolic 3\textendash manifolds $X$. Such an $X$ can be realized as a degree 3 cover of $S^3$ branched over a knot $K$ (\cite{Hi}, \cite{Mo}). \cite{Gr-Gu} show that the distortion $\delta (K)$ of this $K$ satisfies $\delta (K) \gg \text{Vol} (X) \lambda_1 (X)$, (the implied constant being universal). From this and the lower bound for $\lambda_1$ when $X$ varies over congruence arithmetic 3-manifolds (or a congruence thin tower and using the fundamental expansion theorem) one concludes that $K$ and all knots isotopic to it has arbitrarily large distortion by choosing such $X$ of large volume.
 \subsection{Rotation Groups}
Let $\varGamma = \langle \sigma_1 , \sigma_2 , \cdots , \sigma_t \rangle$ be a finitely generated subgroup of the group $SO_3 (\mathbb{R} )$. There is an archimidedian analogue of the expander property for the congruence graphs in this setting and which likewise has many applications (\cite{Lu3}, \cite{Sa3}). Define $T_\sigma$ to be the averaging operator on functions on the two sphere $S^2$ by:
\begin{equation}
   T_\sigma f (x) = \sum^t_{j=1}\limits [f(\sigma_j x) + f(\sigma_j^{-1} x)] .
\end{equation}
$T_\sigma$ is self adjoint on $L^2 (S^2 , dA)$, where $dA$ is the rotation invariant area element on $S^2$, and its spectrum is contained in $[-2t , 2t]$. The spectral gap property is that $2t$ (which is an eigenvalue with eigenvector the constant function) is a simple and isolated point of the spectrum. It is not hard to see that this property depends only on $\varGamma$ and not on the generators. It is conjectured that $\varGamma$ has such a spectral gap iff the $Zcl(\varGamma ) = SO_3$ (which in this case is equivalent to the topological closure of $\varGamma$ being $SO_3 (\mathbb{R} )$). A lot is known towards this conjecture. The first example of a $\varGamma$ with a spectral gap was given by Drinfeld \cite{Dr} and this provided the final step in the solution of the Ruziewicz problem; that the only finitely additive rotationally invariant measure defined on Lebesque measurable subsets of $S^2$, is a multiple of $dA$. His proof of the spectral gap makes use of an arithmetic such $\varGamma$ together with the full force of automorphic forms and the solution of the Ramanujan Conjectures for holomorphic cusp forms on the upper half plane. In \cite{G-J-S} many thin $\varGamma$'s are shown to have a spectral gap. The best result known is the analogue of the Fundamental Expansion Theorem in this context (\cite{B-G2}, \cite{B-G3}), and it suffices for most applications. It asserts that if the matrix elements of members of $\varGamma$ are algebraic, then the Conjecture is true for $\varGamma$. Like the very thin cases of the Fundamental Expansion Theorem, part of the proof here relies on additive combinatorics. This time one needs the full force of the proof of the local Erdos-Volkmann ring conjecture (\cite{E-M}, \cite{Bo})\textemdash that a subset of $\mathbb{R}$ which is closed under addition and multiplication has Hausdorff dimension zero or one. As far as some concrete applications of the spectral gap for these groups, we mention the speed of equidistribution of directions associated with general quaquaversal tilings of three dimensional space (\cite{D-S-V}, \cite{R-S1}) and constructions of quantum gates in the theory of quantum computation (the Solovay-Kitaev Theorem; see \cite{H-R-C}).
\section{Ubiquity of Thin Groups}
Given a finitely generated group $\varGamma$ in $GL_n (\mathbb{Z} )$, one can usually compute $G = Zcl (\varGamma )$ without too much difficulty. On the other hand, deciding if $\varGamma$ is thin can be formidable. In fact one is flirting here with questions that have no decision procedures (I thank Rivin for alerting me to these pitfalls that are close by). For example if \linebreak
$\varGamma = SL_2 (\mathbb{Z} ) \times  SL_2 (\mathbb{Z})$ then there is no decision procedure to determine if an element $A \in \varGamma$ is in the group generated by a general set of say seven elements \cite{Mi}. Even for Gromov hyperbolic groups, the question of whether a finitely generated subgroup generates a finite index subgroup, has no decision procedure (\cite{Rip}, \cite{B-M-S}). Mercifully strong and superstrong approximation only ask about $Zcl (\varGamma )$. Still one is curious about thinness when applying these theorems and sometimes for good reason. For example, in the affine sieve setting, the quality of the expansion impacts the results dramatically (see \cite{N-S} for the cases when $\varGamma$ is a lattice) while the diophantine orbit problems become more standard ones of integer points on homogeneous varieties, when $\varGamma$ is a lattice. Whether the typical $\varGamma$ is thin or not is not so clear, and may depend on how $\varGamma$ arises.
\subsection{Schottky, Ping-Pong:}
Schottky groups in which the generators play ping-pong (\cite{T}, \cite{B-Ge}) are one of the few classes of discrete groups whose group theoretic structure is very simple. If one chooses $A_1, A_2, \ldots , A_\ell$ independently and at random in $SL_n (\mathbb{Z}) (n \geq 2)$, then with high probability $\varGamma = \langle A_1 , \ldots , A_\ell \rangle$ will be free on these generators, Zariski dense in $SL_n$ and thin. If the $A_j$'s are chosen at the $m$\textendash th step of a $\mu$\textendash random walk $(m \rightarrow \infty )$ and support $(\mu )$ generates $SL_n (\mathbb{Z} )$, then this was proved in \cite{Ao}. A more geometric version is proven in \cite{F-R} where the $A$'s are chosen independently and uniformly by taking them from the set of $B$'s with $\max ( \| B \| , \ \| B^{-1} \| )$ less than $X$. Here $\| \ \|$ is any Euclidian norm on the space of matrices and $X \longrightarrow \infty$. Not only is $\varGamma$ thin but it is very thin in the sense that the Hausdorff dimension of the limit set of $\varGamma$ acting on $\mathbb{P}^{n-1} (\mathbb{R})$ is arbitrarily small.

\subsection{Nonarithmetic Lattices}
If $\varGamma \leq G$ \ with $G \neq SL_2 (\mathbb{R})$, is an irreducible nonarithmetic lattice in a semisimple real group $G$, then $\varGamma$ is naturally thin in the appropriate product by its conjugates. The argument is the same as the one in Section 2.3 using local rigidity. The certificate of being thin is that $\varGamma$ is discrete in the factor corresponding to $G$. Examples of this kind which come from monodromy of hypergeometric differential equations in several variables are given in \cite{D-M} and in one variable in \cite{C-W}. It appears that these were the first examples of thin monodromy groups (see section 3.5 below). Other examples of thin monodromy groups in products of $SL_2$'s are given in \cite{N2} and these examples aren't even finitely presented. Teichmuller curves in the moduli space $M_2$ of curves of genus 2, give via Abel-Jacobi, curves in $A_2$ whose monodromies (inclusion of fundamental groups) are thin (\cite{Mc1}). Here too the Zariski closure is a nontrivial product in $Sp_4$ and the thinness follows from having a discrete projection.

\subsection{Reflection Groups in Hyperbolic Space}
Let $f$ be an integral quadratic form in $n$-variables and of signature $(n-1, 1)$. For $n \geq 3, \ O_f (\mathbb{Z})$ the group of integral automorphs of $f$ is a lattice in $G = O_f (\mathbb{R})$. The reflective subgroup $R_f$ is the subgroup of $O_f (\mathbb{Z})$ which is generated by all the hyperbolic reflections which are in $O_f (\mathbb{Z})$. $R_f$ is a normal subgroup of $O_f (\mathbb{Z})$ and if it is nontrivial, then $Zcl(R_f ) = O_f$. Vinberg \cite{Vi} and Nikulin \cite{Ni} have examined the question of when $R_f$ is of finite index in $O_f (\mathbb{Z})$ (they call such an $f$ reflective). In particular, in \cite{Ni} it is shown that there are only finitely many $f$'s (up to integral equivalence) which are reflective. Thus for all but finitely many $f$'s, $R_f$, if it is nontrivial, is a thin group in $GL_n (\mathbb{Z})$ (albeit infinitely generated). Note that Nikulin's theorem fails for $n=2$. If $f$ is a binary quadratic form, then $f$ is reflective iff it is ambiguous in the sense of Gauss (see \cite{Sa4}) and Gauss determined the ambiguous forms in his study of genus theory.

\subsection{Rotation Groups}
An interesting family of rotation groups are the groups \linebreak
$\varGamma (m, n), \ m \geq 3, \ \ n \geq 3$ generated by $\sigma_m$ and $\tau_n$ where
\begin{equation*}
 \sigma_m =
\begin{bmatrix}
   \cos 2\pi /m &\sin 2\pi /m &0\\
- \sin 2\pi / m & \cos 2\pi / m &0\\
0 &0 &1
\end{bmatrix} , \ \ \tau_n = 
\begin{bmatrix}
   1 &0 &0\\
0 &\cos 2\pi / n &\sin 2\pi / n\\
0 &-\sin 2\pi / n &\cos 2\pi / n
\end{bmatrix} .
\end{equation*}
 That is $\varGamma (m,n)$ is a subgroup of $SO_f (\mathbb{R}), \ f(x_1,  x_2 ,  x_3 ) = x^2_1 + x^2_2 + x^2_3$, generated by two rotations about orthogonal axes and of orders $m$ and $n$ respectively. These arise in the theory of quaquaversal tilings of 3\textendash space and their generalizations (\cite{C-R}, \cite{R-S1}).

As abstract groups, these are free products of two cyclic (or dihedral) groups amalgamated over a similar such group (except for $\varGamma (4,4)$ which is finite and which we avoid), see \cite{R-S2}. This description can be used to decide the question of whether $\varGamma (m,n)$ is thin or not and also to show that thin is the rule rather than the exception. If $K = \mathbb{Q} (\cos 2\pi / m,   \sin 2\pi / m,  \cos 2\pi / n,  \sin 2\pi / n)$, then $K$ is a totally real Galois extension of $\mathbb{Q}$ with abelian Galois group $G_{m,n}$. It is plain that $\varGamma (m,n)$ is a subgroup of $SO_f (O[\frac{1}{2}])$, where $O$ is the ring of integers of $K$. Moreover, since $\varGamma(m,n)$ is infinite, the powers of 2 in the denominators of the matrix entries of $\varGamma(m,n)$ must be unbounded (otherwise $\varGamma (m,n)$ would be a discrete subgroup of the compact group $\displaystyle \prod_{\upsilon | \infty}\limits SO_f (K_v )$). Hence the smallest $S$-arithmetic group to contain a subgroup commensurable with $\varGamma (m,n)$ is $SO_f (O_S )$ where $O_S$ are the $S$-integers of $K$, and $S$ consists of the places of $K$ dividing 2. Our thinness question is whether $\varGamma (m,n)$ is of finite or infinite index in the latter. If $| S | \geq 2$, then any finite index subgroup of $SO_f (O_S )$ is a lattice in the higher rank group, $\displaystyle \prod_{\upsilon | (2)}\limits SO_f (K_v )$. By well known rigidity properties of such lattices \cite{Ma} (or one can argue with vanishing of first cohomology groups) and the description of $\varGamma (m,n)$ mentioned above, it follows that $\varGamma (m,n)$ cannot be such a lattice. That is if $| S | \geq 2$, then $\varGamma (m,n)$ is thin and the former holds most of the time (for example if $G_{m,n}$ is not cyclic then $| S | \geq 2)$. If $| S | = 1$, then $\varGamma (m,n)$ may be arithmetic and it is so in some special cases.\footnote{The quaquaversal tiling \cite{C-R} has symmetry group $\varGamma (3,6)$ which is arithmetic \cite{Ser}, while the Dite/Kart tiling \cite{R-S1} has symmetry $\varGamma (10,4)$ for which $K = \mathbb{Q} (\cos \frac{\pi}{10}), G_{10,4} = \mathbb{Z} / 2\mathbb{Z} \times \mathbb{Z} /2 \mathbb{Z}$ and $| S | = 2$, hence the latter is thin.} Perhaps the most interesting cases where $| S | = 1$ are when $m=4$ and $n=2^\nu , \ \nu \geq 3$, \ for which 2 is totally ramified. These have been investigated in \cite{Ro} and \cite{Ser}. Serre shows that for $\nu = 3$ and 4, \ $\varGamma (4, 2^\nu )$ is arithmetic (in fact $\varGamma (4, 2^\nu ) = SO_f (O [\frac{1}{2}])$) while for $\nu \geq 5$, it is thin. The thinness is proven by comparing the Euler characteristics $\chi (SO_f (O [\frac{1}{2}]))$ and $\chi(\varGamma (4,2^\nu ))$, the first using a Tamagawa number computation and the second from the abstract group description of $\varGamma (4,2^\nu )$.

\subsection{Mondromy Groups}
The oldest and perhaps most natural source of finitely generated linear groups comes from monodromy in all of its guises. These include the very classical case of monodromy of the hypergeometric differential equation which we discuss further below, as well as that of a family of varieties varying over a base with its monodromy action on cohomology. For large families, and in cases where the monodromy has been computed, it appears almost always to be arithmetic. The question as to whether such monodromy groups are arithmetic was first raised in \cite{Gr-Sc}. For example for the universal family of smooth projective hypersurfaces of degree $d$ and dimension $n$ in projective space, the monodromy representation on $H^n (X_0 , \mathbb{Z}), X_0$ a base hypersurface, is an arithmetic subgroup of $GL(H^n (X_0 ))(\mathbb{Z})$; see [Be] where the exact level in $G(\mathbb{Z})$ is determined. For smaller families such as cyclic covers of $\mathbb{P}^1$, which have recently been studied in \cite{Mc2} in connection with the thinness question, the story is similar. More precisely, consider the family of curves (in affine coordinates) given by
\begin{equation}
   C_a : y^d = (x - a_1 )(x - a_2 ) \cdots (x - a_{n+1}),  
\end{equation}
where the parameters $a$ vary so that $a_i \neq a_j , \text{ for } i \neq j$. The fundamental group of the space of $a$'s is the pure braid group and it has a monodromy representation on $H_1 (C, \mathbb{Z}) \cong \mathbb{Z}^{2g}, \ g$ the genus of $C_a$, and again $C$ is a fixed base curve. Answering a question in \cite{Mc2}, \cite{Ve1} shows that if $n \geq 2d$, then the image of the monodromy representation of the braid group in $GL(H_1 (C))(\mathbb{Z})$ is arithmetic. This generalizes a result of \cite{A'c} for $d=2$. The proof is based on another result of \cite{Ve2} which asserts that for $\mathbb{Q}$ rank two or higher arithmetic groups, a Zariski dense subgroup which contains enough elements from opposite horospherical  subgroups is necessarily arithmetic. If $n < 2d$, then as observed in \cite{Mc2}, there are examples based on the nonarithmetic lattices of \cite{D-M} in SU(2,1) which are thin (one such is $n=3$ and $d=18)$.\\
\indent The thinness story for monodromy groups of one parameter families is less clear. We discuss in some detail the very rich examples of the classical hypergeometric equation. Let $\alpha , \ \beta \in \mathbb{Q}^n$ and consider the ${_n}F_{n-1}$ algebraic hypergeometric equation:
\begin{equation}
   Du = 0
\end{equation}
where $\displaystyle D = (\theta + \beta_1 -1)(\theta + \beta_2 -1) \cdots (\theta - \beta_n -1) - z (\theta + \alpha_1 ) \cdots (\theta + \alpha_n )$ and $\theta = z \frac{d}{dz}$.\\
\indent The equation is regular outside $\{0, 1, \infty \}$ and the fundamental group $\pi_1 (\mathbb{P}^1 - \{ 0, 1, \infty \} )$ has a representation in $GL_n$ gotten by analytic continuation of a basis of solutions to (3.2) along curves in the thrice punctured sphere. Its image in $GL_n$ is denoted by $H(\alpha , \beta )$ and is the monodromy group in question (defined up to conjugation in $GL_n). \  H(\alpha , \beta )$ is generated by the local monodromies $A, B, C \ \ (C = A^{-1} B)$ gotten from loops about $0, \infty$ and 1 respectively, see Beukers and Heckman \cite{B-H} for a detailed description. We restrict to $H$'s which can be conjugated into $GL_n (\mathbb{Z})$, which is equivalent to the characteristic polynomials of $A$ and $B$ being products of cyclotomic polynomials.\footnote{We assume further that $(\alpha , \beta )$ are primitive in the sense of \cite{B-H}.} Such $H(\alpha , \beta )$'s are self-dual and according to \cite{B-H}, their Zariski closures $G(\alpha , \beta )$, are either finite, $O_n$ or $Sp_n$, and they determine which it is explicitly in terms of $\alpha$ and $\beta$. Our interest is whether $H(\alpha , \beta )$ is of finite or infinite index in $G(\alpha , \beta ) (\mathbb{Z})$. Other than the cases where $H(\alpha , \beta )$ (or equivalently $G(\alpha , \beta )$) are finite, all of which are listed in \cite{B-H}, there are few cases where $H(\alpha , \beta )$ itself is known.\\
\indent Recently Venkatamarana \cite{Ve1} has shown that for $n$ even and 
\begin{align}
   &\alpha = \left( \frac{1}{2} + \frac{1}{n+1} , \frac{1}{2} + \frac{2}{n+1} , \cdots , \frac{1}{2} + \frac{n}{n+1} \right)\notag\\
&\beta = \left( 0, \frac{1}{2} + \frac{1}{n} , \frac{1}{2} + \frac{2}{n} , \cdots , \frac{1}{2} + \frac{n-1}{n} \right) ,
\end{align}
$H(\alpha , \beta )$ is arithmetic (here $G(\alpha , \beta ) = Sp(n)$). He deduces this by showing that for these exact parameters, the monodromy representation of $\pi_1 (\mathbb{P}^1  - \{ 0, 1, \infty \})$ factors through a representation of the braid group on (3.1) with $d=2$. In particular the arithmeticity follows from the arithmeticity of the latter.\\
\indent The very fruitful Dwork family (see \cite{Ka}, \cite{H-S-T}) $n \geq 4$ \ even, and 
\begin{align}
   &\alpha = (0, 0, \cdots , 0) ,\notag \\
&\beta = \left( \frac{1}{n+1} , \frac{2}{n+1}, \cdots , \frac{n}{n+1} \right) ,
\end{align}
is apparently different. Again $G(\alpha , \beta ) = Sp(n)$ and for $n=4$, the local monodromies are 
\begin{equation}
 A =
\begin{bmatrix}
   0 &0 &0 &-1\\
1 &0 &0 &-1\\
0 &1 &0 &-1\\
0 &0 &1 &-1
\end{bmatrix} \ \ \text{ and } \ \ \ C =
\begin{bmatrix}
  1 &0 &0 &5\\
0 &1 &0 &-5\\
0 &0 &1 &5\\
0 &0 &0 &1
\end{bmatrix} .
\end{equation}

Very recently \cite{B-T} have shown that $A$ and $C$ in (3.5) play generalized ping-pong on certain subsets of $\mathbb{P}^3$, from which it follows that $H(\alpha , \beta )  \cong \mathbb{Z} / 5 \mathbb{Z} \ast \mathbb{Z}$. From rigidity, or the first cohomology properties of finite index subgroups of $Sp(4, \mathbb{Z} )$, it follows that $H(\alpha , \beta )$ must be thin. It seems likely that $H(\alpha , \beta )$ is thin for the whole Dwork family, i.e, $n \geq 4$, but other than showing that the corresponding $A$ and $C$'s play ping-pong, there appear to be no known means of proving this and no infinite family of thin $H(\alpha , \beta )$'s with $G(\alpha , \beta )$ symplectic is known. For $n=4$ there are 112 such $H(\alpha , \beta )$'s in $Sp (4,\mathbb{Z} )$ \cite{Si-Ve}. Using extensions of the technique in \cite{Ve1} it is shown in \cite{Si-Ve} that of these, 63 are arithmetic.  Of these 3 [namely $\left( (0, 0, 0, 0), \left( \frac{1}{6}, \frac{1}{6}, \frac{5}{6}, \frac{5}{6}\right)\right), \left( (0,0,0,0),\left(\frac{1}{6}, \frac{1}{4},\frac{3}{4},\frac{5}{6}\right)\right)$ and $\left( (0,0,0,0),\left(\frac{1}{10}, \frac{3}{10}, \frac{7}{10},\frac{9}{10}\right)\right) $] correspond to the 14 hypergeometrics associated with certain Calabi-Yau three folds \cite{C-Y-Y}. Of the other 11, 7 are shown to be thin in \cite{B-T}, again by finding ping-pong sets in $\mathbb{P}^3$. This leaves 4 of these Calabi-Yau's for which the thinness question is open. It would be interesting to understand the geometric significance, if there is one, for $H(\alpha , \beta )$ being thin or not in these families.\\

What is lacking above is a certificate for $H(\alpha , \beta )$ being thin that can be applied for example to families (i.e. $n \rightarrow \infty )$. A robust such certificate has been provided in the case that $G(\alpha , \beta )(\mathbb{R})$ is of rank one and $n > 3$ \cite{F-M-S}.
In these cases $G(\alpha , \beta )$, as a group defined over $\mathbb{Q}$ is $O_f$, where $f$ is a rational quadratic form in an odd number of variables and of signature  $(n-1 , 1)$ (over $\mathbb{R}$). We call these $(\alpha , \beta )$'s hyperbolic hypergeometrics and besides a (long) list of sporadic examples, they come in seven infinite parametric families \cite{F-M-S}. Our conjecture for these is that thin rules, that is for all but finitely many of the hyperbolic hypergeometrics, $H(\alpha , \beta )$ is thin. This is proved in \cite{F-M-S} for a number (but not all) of the seven families. For example for $n$ odd consider the two families:
\begin{align}
   &\alpha = \left( 0, \frac{1}{n+1} , \frac{2}{n+1} , \cdots , \frac{n-1}{2(n+1)} ,  \frac{n+3}{2(n+1)} , \cdots , \frac{n}{n+1} \right)\notag\\
&\beta = \left( \frac{1}{2} , \frac{1}{n} , \frac{2}{n} , \cdots , \frac{n-1}{n} \right),
\end{align}
and 
\begin{align}
   &\alpha = \left( \frac{1}{2}, \frac{1}{2n-2} , \frac{3}{2n - 2} ,   \cdots  \frac{2n-3}{2n-2} \right)\notag\\
&\beta = \left( 0, 0, 0, \frac{1}{n-2} , \frac{2}{n-2} , \cdots \frac{n-3}{n-2} \right) .
\end{align}
Both of these families are hyperbolic hypergeometrics and for both $H(\alpha , \beta )$ is thin for $n \geq 5$ and is arithmetic for $n=3$.\\
\indent The proof is based on the following principle: if $\psi : G(\mathbb{Z} ) \rightarrow K$ is a morphism onto a group $K$ for which $| \psi (H ( \alpha , \beta)) \backslash K | = \infty$, then certainly $H(\alpha , \beta )$ is of infinite index in $G(\mathbb{Z})$. Now in the higher rank cases there are no useful such $\psi$'s (by the Margulis normal subgroup theorem \cite{Ma} in these cases if $K$ is infinite then $ker (\psi )$ is finite), however, in the rank one case such $\psi$'s may exist and yield a certificate of thinness. Indeed in this hyperbolic case if $R_f$ is the Vinberg reflection subgroup described in 3.3, then as mentioned there, except for finitely many $f$'s,  $K_f : = O_f (\mathbb{Z}) / R_f$ is infinite. To use this one needs to analyze the image of $H(\alpha , \beta )$ in $K_f$. The key observation is that up to the finite index the hyperbolic hypergeometrics are generated by Cartan Involutions.\footnote{The local monodromy $C$ about 1 is always a pseudo reflection and in these cases yields a Cartan involution.} These are linear reflections of $\mathbb{Q}^n$ which induce isometries on hyperbolic space given by geodesic inversions in a point [the hyperbolic reflections are generated by root vectors $v$ in $\mathbb{Z}^n$ outside the light cone $(f(v) > 0)$ while the Cartan involutions by root vectors $w$ in $\mathbb{Z}^n$ inside the light cone, in fact $f(w) = -2 ]$. In order to examine the image of a group generated by such Cartan involutions in $K_f$, consider the \textquotedblleft minimum distance graph,\textquotedblright \ $X_f$. Its vertices are the integral Cartan root vectors $V_{-2} (\mathbb{Z}) = \left\{ v \in \mathbb{Z}^n : f(v) = -2 \right\}$, and $v$ and $w$ are joined if $f(v,w) = -3$. One can show that the components of $X_f$ consist of finitely many isomorphism types and each is the Cayley graph of a finitely generated Coxeter group. The main lemma \cite{F-M-S} asserts that if $\Sigma \subset V_{-2} (\mathbb{Z})$ is a connected component of $X_f$ then the image of the group generated by the Cartan involutions $r_v$ with roots $v \in \Sigma$, is a finite subgroup of $K_f$.\footnote{The proof makes use of the quite special feature of the binary form $g = x^2 + 3 xy + y^2$, of being integrally equivalent to $-g$ (called reciprocal in \cite{Sa4}).} This together with Vinberg and Nikulin's theorems gives a robust certificate for the thinness of these hyperbolic hypergeometric monodromies. As far as I know (3.6) and (3.7) give the first family of thin monodromy groups in high dimensions for which $G$ is simple.\\
\indent We end with some comments about the arithmetic Ramanujan Conjectures. The gonality of a congruence arithmetic surface being linear in its genus and the Heegaard genus of a congruence hyperbolic three manifold being linear in its volume, as well as the proof that there are only finitely many maximal arithmetic reflection groups, all appeal to the uniform lower bounds for $\lambda_1$ for all such manifolds. This follows from what is known towards the Ramanujan Conjectures but it does not follow from the fundamental expansion theorem since the latter only applies to one tower at a time. As far as the general Ramanujan Conjectures, some progress has been made since the report \cite{Sa1}. Namely in \cite{Ar} a precise formulation of the Ramanujan Conjectures for these groups is given, and moreover it is shown (assuming forms of the fundamental lemma which themselves should be theorems before too long) that these conjectures will follow if one can prove the Ramanujan Conjectures for $GL_m$. \\

\noindent \underline{Acknowledgement}: These brief notes cover a lot of ground. I thank my collaborators, the people whose work is quoted and the many mathematicians with whom I have discussed aspects of the theory connected with these thin groups. Thanks to the referee for pointing me to the relevant applications in \cite{Gr-Gu}.

\end{document}